\documentclass[11pt,bezier]{article}
\usepackage{amsmath}
\usepackage{amsfonts,amsthm,amssymb}
\usepackage{amsfonts}
\usepackage{graphics}
\usepackage{amssymb}
\newtheorem{lem}{Lemma}[section]
\newtheorem{thm}[lem]{Theorem}
\newtheorem{cor}[lem]{Corollary}
\newtheorem{conj}{Conjecture}

\theoremstyle{definition}

\begin{document}
\title{On the sizes of vertex-$k$-maximal $r$-uniform hypergraphs
\footnote{The research is supported by NSFC (Nos. 11531011, 11771039, 11771443).}}
\author{Yingzhi Tian$^{a}$ \footnote{Corresponding author. E-mail: tianyzhxj@163.com (Y. Tian), hjlai@math.wvu.edu (H. Lai), mjx@xju.edu.cn (J. Meng).}, Hong-Jian Lai$^{b}$, Jixiang Meng$^{a}$ \\
{\small $^{a}$College of Mathematics and System Sciences, Xinjiang
University, Urumqi, Xinjiang 830046, PR China}\\
{\small $^{b}$Department of Mathematics, West Virginia University,
Morgantown, WV 26506, USA}}

\date{}

\maketitle

\noindent{\bf Abstract } Let $H=(V,E)$ be a hypergraph, where $V$ is a  set of vertices and $E$ is a set of non-empty subsets of $V$ called edges. If all edges of $H$ have the same cardinality $r$, then $H$ is a $r$-uniform hypergraph; if $E$ consists of all $r$-subsets of $V$, then $H$ is a complete $r$-uniform hypergraph, denoted by $K_n^r$, where $n=|V|$. A hypergraph $H'=(V',E')$ is called a subhypergraph of $H=(V,E)$ if $V'\subseteq V$ and $E'\subseteq E$. A $r$-uniform hypergraph $H=(V,E)$ is vertex-$k$-maximal if every subhypergraph of $H$ has vertex-connectivity at most $k$, but for any edge $e\in E(K_n^r)\setminus E(H)$, $H+e$ contains at least one subhypergraph with vertex-connectivity at least $k+1$. In this paper, we first prove that for given integers $n,k,r$ with $k,r\geq2$ and $n\geq k+1$,  every vertex-$k$-maximal $r$-uniform hypergraph $H$ of order $n$ satisfies $|E(H)|\geq (^n_r)-(^{n-k}_r)$, and this lower bound is best possible. Next, we conjecture that for sufficiently large $n$, every vertex-$k$-maximal $r$-uniform hypergraph $H$ on $n$ vertices satisfies $|E(H)|\leq(^n_r)-(^{n-k}_r)+(\frac{n}{k}-2)(^k_r)$, where $k,r\geq2$ are integers. And the conjecture is verified for the case $r>k$.


\section{Introduction}

In this paper, we consider finite simple graphs. For graph-theoretical terminologies and notation not defined here, we follow \cite{Bondy}.
For a graph $G$, we use $\kappa(G)$ to denote the $vertex$-$connectivity$ of $G$. The $complement$ of a graph $G$ is denoted by $G^c$. For $X\subseteq E(G^c)$, $G+X$ is the graph with vertex set $V(G)$ and edge set $E(G)\cup X$. We will use $G+e$ for $G+\{e\}$. The $floor$ of a real number $x$, denoted by $\lfloor x\rfloor$, is the greatest integer not larger than $x$; the $ceil$ of a real number $x$, denoted by $\lceil x\rceil$, is the least integer greater than or equal to $x$. For two integers $n$ and $k$, we define $(_k^n)=\frac{n!}{k!(n-k)!}$ when $k\leq n$ and $(_k^n)=0$ when $k>n$.

Matula \cite{Matula78} first explicitly studied the quantity
$\overline{\kappa}(G)=max\{\kappa(G'): G'\subseteq G\}$. For a positive integer $k$, the graph $G$ is $vertex$-$k$-$maximal$ if $\overline{\kappa}(G)\leq k$ but for any edge $e\in E(G^c)$, $\overline{\kappa}(G+e)>k$. Because $\kappa(K_n)=n-1$, a vertex-$k$-maximal graph $G$ with at most $k+1$ vertices must be a complete graph.

The $union$ of two graphs $G_1$ and $G_2$, denoted by $G_1\cup G_2$, is the graph with vertex set $V(G_1)\cup V(G_2)$ and edge set $E(G_1)\cup E(G_2)$.
The $join$ of two graphs $G_1$ and $G_2$, denoted by $G_1\vee G_2$, is the graph obtained from the union of $G_1$ and $G_2$ by adding all the edges that connect the vertices of $G_1$ with $G_2$.
Let $G_{n,k}=((p-1)K_k\cup K_q) \vee K_k^c$, where $n=pk+q\geq2k$ ($1\leq q\leq k$) and $(p-1)K_k$ is the union of $p-1$ complete graphs on $k$ vertices.  Then $G_{n,k}$ is vertex-$k$-maximal and  $|E(G_{n,k})|\leq\frac{3}{2}(k-\frac{1}{3})(n-k)$, where the equality holds if $n$ is a multiple of $k$.
Mader \cite{Mader79} conjectured that, for large order of graphs, the graph $G_{n,k}$ would in fact present the best possible upper bound for the sizes of a vertex-$k$-maximal graph.

\begin{conj} (Mader \cite{Mader79})
Let $k\geq2$ be an integer. Then for sufficiently large $n$, every vertex-$k$-maximal graph on $n$ vertices satisfies $|E(G)|\leq\frac{3}{2}(k-\frac{1}{3})(n-k)$.
\end{conj}

Some progresses towards Conjecture 1 are listed in the following.

\begin{thm}
Let $k\geq2$ be an integer.

($i$) (Mader \cite{Mader72}, see also \cite{Mader79}) Conjecture 1 holds for $k\leq6$.

($ii$) (Mader \cite{Mader72}, see also \cite{Mader79}) For sufficiently large $n$, every vertex-$k$-maximal graph $G$ on $n$ vertices satisfies $|E(G)|\leq(1+\frac{1}{\sqrt{2}})k(n-k)$.

($iii$) (Yuster \cite{Yuster}) If $n\geq\frac{9k}{4}$, then every vertex-$k$-maximal graph $G$ on $n$ vertices satisfies $|E(G)|\leq\frac{193}{120}k(n-k)$.

($iv$) (Bernshteyn and Kostochka \cite{Bernshteyn}) If $n\geq\frac{5k}{2}$, then every vertex-$k$-maximal graph $G$ on $n$ vertices satisfies $|E(G)|\leq\frac{19}{12}k(n-k)$.
\end{thm}

In \cite{Xu}, Xu, Lai and Tian obtained the lower bound of the sizes of  vertex-$k$-maximal graphs.

\begin{thm} (Xu, Lai and Tian \cite{Xu})
Let $n,k$ be integers with $n\geq k+1\geq3$. If $G$ is a vertex-$k$-maximal graph on $n$ vertices, then $|E(G)|\geq (n-k)k+\frac{k(k-1)}{2}$. Furthermore, this bound is best possible.
\end{thm}

The related studies on edge-$k$-maximal graphs have been conducted by quite a few researchers, as seen in [7,9,12,13,15], among others. For corresponding digraph problems, see [1,8], among others.

Let $H=(V,E)$ be a hypergraph, where $V$ is a finite set and $E$ is a set of non-empty subsets of $V$, called edges. An edge of cardinality 2 is just a graph edge. For a vertex $u\in V$ and an edge $e\in E$, we say $u$ is $incident$ $with$ $e$ or $e$ is $incident$ $with$ $u$ if $u\in e$.
If all edges of $H$ have the same cardinality $r$, then $H$ is a $r$-$uniform$ $hypergraph$; if $E$ consists of all $r$-subsets of $V$, then $H$ is a $complete$ $r$-$uniform$ $hypergraph$, denoted by $K_n^r$, where $n=|V|$.
For $n<r$, the complete $r$-uniform hypergraph $K_n^r$ is just the hypergraph with $n$ vertices and no edges.
The $complement$ of a $r$-uniform hypergraph $H=(V,E)$, denoted by $H^c$, is the $r$-uniform hypergraph with vertex set $V$ and edge set consisting of all $r$-subsets of $V$ not in $E$. A hypergraph $H'=(V',E')$ is called a $subhypergraph$ of $H=(V,E)$, denoted by $H'\subseteq H$, if $V'\subseteq V$ and $E'\subseteq E$.
For $X\subseteq E(H^c)$, $H+X$ is the hypergraph with vertex set $V(H)$ and edge set $E(H)\cup X$; for $X'\subseteq E(H)$, $H-X'$ is the hypergraph with vertex set $V(H)$ and edge set $E(H)\setminus X'$. We use $H+e$ for $H+\{e\}$ and $H-e'$ for $H-\{e'\}$ when $e\in E(H^c)$ and $e'\in E(H)$.
For $Y\subseteq V(H)$, we use $H[Y]$ to denote the hypergraph $induced$ by $Y$, where $V(H[Y])=Y$ and $E(H[Y])=\{e\in E(H): e\subseteq Y\}$. $H-Y$ is the hypergraph induced by $V(H)\setminus Y$.

Let $H$ be a hypergraph and $V_1,V_2,\cdots, V_l$ be subsets of $V(H)$. An edge $e\in E(H)$ is $(V_1,V_2,\cdots, V_l)$-$crossing$ if $e\cap V_i\neq\emptyset$ for $1\leq i\leq l$. If in addition, $e\subseteq \cup_{i=1}^lV_i$, then $e$ is $exact$-$(V_1,V_2,\cdots, V_l)$-$crossing$. The set of all $(V_1,V_2,\cdots, V_l)$-crossing edges of $H$ is denoted by $E_H[V_1,V_2,\cdots, V_l]$; the set of all exact-$(V_1,V_2,\cdots, V_l)$-crossing edges of $H$ is denoted by $E_{H[V_1\cup V_2\cup\cdots\cup V_l]}[V_1,V_2,\cdots, V_l]$. Let $d_H(V_1,V_2,\cdots, V_l)=|E_H[V_1,V_2,\cdots, V_l]|$ and $d_{H[V_1\cup V_2\cup\cdots\cup V_l]} \\ (V_1,V_2,\cdots, V_l)=|E_{H[V_1\cup V_2\cup\cdots\cup V_l]}[V_1,V_2,\cdots, V_l]|$. For a vertex $u\in V(H)$, we call $d_H(u):=d_H(\{u\}, V(H)\setminus \{u\})$ the $degree$ of $u$ in $H$.  The $minimum$ $degree$ $\delta(H)$ of $H$ is defined as $min\{d_H(u): u\in V\}$; the $maximum$ $degree$ $\Delta(H)$ of $H$ is defined as $max\{d_H(u): u\in V\}$. When $\delta(H)=\Delta(H)=k$, we call $H$ $k$-$regular$.

Given a hypergraph $H$, we define a $walk$ in $H$ to be an  alternating sequence $v_1,e_1,v_2,\cdots,e_s,\\ v_{s+1}$ of vertices and edges of $H$ such that: $v_i\in V(H)$ for $i=1,\cdots,s+1$; $e_i\in E(H)$ for $i=1,\cdots,s$; and $v_i,v_{i+1}\in e_i$ for $i=1,\cdots,s$. A $path$ is a walk with additional restrictions that the vertices are all distinct and the edges are all distinct. A hypergraph $H$ is $connected$ if for every pair of vertices $u,v\in V(H)$, there is a path connecting $u$ and $v$; otherwise $H$ is $disconnected$. A $component$ of a hypergraph $H$ is a maximal connected subhypergraph of $H$. A subset $X\subseteq V$ is called a $vertex$-$cut$ of $H$ if $H-X$ is disconnected. We define the $vertex$-$connectivity$ of $H$, denoted by $\kappa(H)$, as follows: if $H$ had at least one vertex-cut, then
$\kappa(H)$ is the cardinality of a minimum vertex-cut of $H$; otherwise $\kappa(H)=|V(H)|-1$.  We call a hypergraph $H$ $k$-$vertex$-$connected$ if $\kappa(H)\geq k$. Let $\overline{\kappa}(H)=max\{\kappa(H'): H'\subseteq H\}$. For a positive integer $k$, the $r$-uniform hypergraph $H$ is $vertex$-$k$-$maximal$ if $\overline{\kappa}(H)\leq k$ but for any edge $e\in E(H^c)$, $\overline{\kappa}(H+e)>k$. Since $\kappa(K_n^r)=n-r+1$, we note that $H$ is  complete if $H$ is a vertex-$k$-maximal $r$-uniform hypergraph with $n-r+1\leq k$, where $n=|V(H)|$. The edge-$k$-maximal hypergraph
can be defined similarly.
For results on the connectivity of hypergraphs, see [2,5,6] for references.

In \cite{Tian}, we determined, for given integers $n$, $k$ and $r$, the extremal sizes of an edge-$k$-maximal $r$-uniform hypergraph on $n$ vertices.

\begin{thm}  (Tian, Xu, Lai and Meng \cite{Tian})
Let $k$ and $r$ be integers with $k,r\geq2$, and let $t=t(k,r)$ be the largest integer such that $(^{t-1}_{r-1})\leq k$. That is, $t$ is the integer satisfying $(^{t-1}_{r-1})\leq k<(^{t}_{r-1})$.
If $H$ is an edge-$k$-maximal $r$-uniform hypergraph with $n=|V(H)|\geq t$, then

($i$) $|E(H)|\leq (^{t}_{r})+(n-t)k$, and this bound is best possible;

($ii$) $|E(H)|\geq (n-1)k -((t-1)k-(^{t}_{r}))\lfloor\frac{n}{t}\rfloor$, and this bound is best possible.
\end{thm}

The main goal of this research is to investigate, for given integers $n$, $k$ and $r$, the extremal sizes of a vertex-$k$-maximal $r$-uniform hypergraph on $n$ vertices.
Section 2 below is devoted to the study of some properties of vertex-$k$-maximal $r$-uniform hypergraphs. In Section 3, we give the best possible lower bound of the sizes of vertex-$k$-maximal $r$-uniform hypergraphs. We propose a conjecture on the upper bound of the sizes of vertex-$k$-maximal $r$-uniform hypergraphs and verify the conjecture for the case $r>k$ in Section 4.

\section{Properties of vertex-$k$-maximal $r$-uniform hypergraphs}

Combining the definition of vertex-$k$-maximal $r$-uniform hypergraph with $\kappa(K_n^r)=n-r+1$, we obtain that $H$ is isomorphic to $K_n^r$ if
$H$ is a vertex-$k$-maximal $r$-uniform hypergraph with $n=|V(H)|\leq k+r-1$.

\begin{lem}
Let $n,k,r$ be integers with $k,r\geq2$ and $n\geq k+r-1$. If
$H$ is a vertex-$k$-maximal $r$-uniform hypergraph on $n$ vertices, then $\overline{\kappa}(H)=\kappa(H)=k$.
\end{lem}

\noindent{\bf Proof.} Since $H$ is vertex-$k$-maximal, we have $\kappa(H)\leq\overline{\kappa}(H)\leq k$. In order to complete the proof, we only need to show that $\kappa(H)\geq k$.

If $n=k+r-1$, then $H$ is complete and  $\kappa(H)=n-r+1=k$. Thus, assume $n\geq k+r$, and so $H$ is not complete. On the contrary, assume $\kappa(H)< k$. Since $H$ is not complete, $H$ has a vertex-cut $S$ with $|S|=\kappa(H)< k$. Let $C_1$ be a component of $H-S$ and $C_2=H-(S\cup V(C_1))$.
By $|V(C_1)\cup V(C_2)|=n-|S|\geq k+r-(k-1)=r+1$, we can choose a $r$-subset $e\subseteq V(C_1)\cup V(V_2)$ such that $e\cap V(C_i)\neq\emptyset$ for $i=1,2$. Then $e\in E(H^c)$.

Since $H$ is vertex-$k$-maximal, we have $\overline{\kappa}(H+e)\geq k+1$. Hence $H+e$ contains a subhypergraph $H'$ with  $\kappa(H')=\overline{\kappa}(H+e)\geq k+1$.  Since $\overline{\kappa}(H)\leq k$, $H'$ cannot be a subhypergraph of $H$, and so $e\in E(H')$. Since $V(H')\cap V(C_i)\neq\emptyset$ for $i=1,2$, it follows that $V(H')\cap S$ is a vertex-cut of $H'-e$.

Since $|V(C_1)\cup V(C_2)|=n-|S|\geq k+r-(k-1)=r+1\geq3$, one of $C_i$, say  $C_1$, contains at least two vertices. Let $u_1\in e\cap V(C_1)$. Then $S'= (V(H')\cap S)\cup\{u_1\}$ is a vertex-cut of $H'$, and so we obtain
$$k+1>|S|+1\geq|V(H')\cap S|+1=|S'|\geq\kappa(H')\geq k+1,$$
a contradiction.
$\Box$

Let $H$ be a vertex-$k$-maximal $r$-uniform hypergraph with $|V(H)|\geq k+r$. By Lemma 2.1, $\overline{\kappa}(H)=\kappa(H)=k$. By $|V(H)|\geq k+r$, $H$ is not complete, thus $H$ contains vertex-cuts.   Let $S$ be a minimum vertex-cut of $H$, $C_1$ be a component of $H-S$ and $C_2=H-(S\cup V(C_1))$. We call $(S,H_1,H_2)$ a $separation$ $triple$ of $H$, where $H_1=H[S\cup V(C_1)]$ and $H_2=H[S\cup V(C_2)]$.

\begin{lem}
Let $n,k,r$ be integers with $k,r\geq2$ and $n\geq k+r$, and
$H$ be a vertex-$k$-maximal $r$-uniform hypergraph on $n$ vertices.
Assume $(S,H_1,H_2)$ is a separation triple  of $H$. If $e\in E(H_1^c)\cup E(H_2^c)$, then any subhypergraph $H'$ of $H+e$ with $\kappa(H')\geq k+1$ is either a subhypergraph of $H_1+e$ or a subhypergraph of $H_2+e$.
Furthermore, if $e\subseteq E(H_i^c)\setminus E((H[S])^c)$, then $H'$ is a subhypergraph of $H_i+e$ for $i=1,2$.
\end{lem}

\noindent{\bf Proof.} Let $e\in E(H_1^c)\cup E(H_2^c)$. Since $H$ is vertex-$k$-maximal, we have $\overline{\kappa}(H+e)\geq k+1$. Let $H'$ be a subhypergraph of $H+e$ with $\kappa(H')=\overline{\kappa}(H+e)\geq k+1$. We assume, on the contrary, that $V(H')\cap (V(H_1)-S)\neq\emptyset$ and $V(H')\cap (V(H_2)-S)\neq\emptyset$. This, together with $e\in E(H_1^c)\cup E(H_2^c)$, implies that $S\cap V(H')$ is a vertex-cut of $H'$. Hence $k=|S|\geq |S\cap V(H')|\geq\kappa(H')\geq k+1$, a contradiction.  Therefore, we cannot have both $V(H')\cap (V(H_1)-S)\neq\emptyset$ and $V(H')\cap (V(H_2)-S)\neq\emptyset$. If $V(H')\cap (V(H_1)-S)=\emptyset$, then $H'$ is a subhypergraph of $H_2+e$; if $V(H')\cap (V(H_2)-S)=\emptyset$, then $H'$ is a subhypergraph of $H_1+e$.

If $e\subseteq E(H_1^c)\setminus E((H[S])^c)$, then $V(H')\cap (V(H_1)-S)\neq\emptyset$ and $V(H')\cap (V(H_2)-S)=\emptyset$, thus $H'$ is a subhypergraph of $H_1+e$. Similarly, if  $e\subseteq E(H_2^c)\setminus E((H[S])^c)$, then  $H'$ is a subhypergraph of $H_2+e$.
$\Box$

\begin{lem}
Let $n,k,r$ be integers with $k,r\geq2$ and $n\geq k+r$, and
$H$ be a vertex-$k$-maximal $r$-uniform hypergraph on $n$ vertices.
Assume $(S,H_1,H_2)$ is a separation triple of $H$ and $n_i=|V(H_i)|$ for $i=1,2$. Then

($i$) $E_{H^c}[V(H_1)-S, S, V(H_2)-S]=\emptyset$, and

($ii$) $d_H(V(H_1)-S, S, V(H_2)-S)=(^{n}_r)-(^{n_1}_r)-(^{n_2}_r)+(^{k}_r)
-(^{n-k}_r)+(^{n_1-k}_r)+(^{n_2-k}_r)$.
\end{lem}

\noindent{\bf Proof.}
($i$) By contradiction, assume $E_{H^c}[V(H_1)-S, S, V(H_2)-S]\neq\emptyset$. Let $e\in E_{H^c}[V(H_1)-S, S, V(H_2)-S]$. Since $H$ is vertex-$k$-maximal, there is a subhypergraph $H'$ of $H+e$ such that $\kappa(H')=\overline{\kappa}(H+e)\geq k+1$. By $\overline{\kappa}(H)\leq k$, $e\in E(H')$. This, together with $e\in E_{H^c}[V(H_1)-S, S, V(H_2)-S]$, implies $V(H')\cap S \neq\emptyset$ and $V(H')\cap(V(H_i)-S)\neq\emptyset$ for $i=1,2$. Hence $S\cap V(H')$ is a vertex-cut of $H'$. But then we obtain $k=|S|\geq |S\cap V(H')|\geq\kappa(H')\geq k+1$, a contradiction. It follows $E_{H^c}[V(H_1)-S, S, V(H_2)-S]=\emptyset$.

($ii$) By ($i$),  $E_{H^c}[V(H_1)-S, S, V(H_2)-S]=\emptyset$. This implies that if $e$ is a $r$-subset such that $e\cap S\neq\emptyset$ and $e\cap (V(H_i)-S)\neq\emptyset$ for $i=1,2$, then $e\in E(H)$. Since the number of $r$-subsets contained in $V(H_1)$ or $V(H_2)$ is $(^{n_1}_r)+(^{n_2}_r)-(^{k}_r)$, and the number of $r$-subsets exactly intersecting $V(H_1)-S$ and $V(H_1)-S$ is $(^{n-k}_r)-(^{n_1-k}_r)-(^{n_2-k}_r)$, we have

\ \ \ \ \ \ \ \ \ \ \ \ \ \ $d_H(V(H_1)-S, S, V(H_2)-S)$

\ \ \ \ \ \ \ \ \ \ \ \ \ \ $=|E_H[V(H_1)-S, S, V(H_2)-S]|$

\ \ \ \ \ \ \ \ \ \ \ \ \ \ $=(^n_r)-((^{n_1}_r)+(^{n_2}_r)-(^{k}_r))-
((^{n-k}_r)-(^{n_1-k}_r)-(^{n_2-k}_r))$

\ \ \ \ \ \ \ \ \ \ \ \ \ \ $=(^{n}_r)-(^{n_1}_r)-(^{n_2}_r)+(^{k}_r)
-(^{n-k}_r)+(^{n_1-k}_r)+(^{n_2-k}_r)$.

This completes the proof.
$\Box$

\section{The lower bound of the sizes of vertex-$k$-maximal $r$-uniform hypergraphs}

The $union$ of two hypergraphs $H_1$ and $H_2$, denoted by $H_1\cup H_2$, is the hypergraph with vertex set $V(H_1)\cup V(H_2)$ and edge set $E(H_1)\cup E(H_2)$.
The $r$-$join$ of two hypergraphs $H_1$ and $H_2$, denoted by $H_1\vee_r H_2$, is the  hypergraph obtained from the union of $H_1$ and $H_2$ by adding all the edges with cardinality $r$ that connect the vertices of $H_1$ with $H_2$.

\noindent{\bf Definition 1.} Let $n,k,r$ be integers such that $k,r\geq2$ and $n\geq k+1$. We define $H_L(n;k,r)$ to be $K_k^r\vee_r (K_{n-k}^r)^c$.

\begin{lem}
Let $n,k,r$ be integers such that $k,r\geq2$ and $n\geq k+1$. If $H= H_L(n;k,r)$, then

($i$) $H$ is vertex-$k$-maximal, and

($ii$) $|E(H)|=(^n_r)-(^{n-k}_r)$.
\end{lem}

\noindent{\bf Proof.} ($i$) By Definition 1, $H$ is obtained from the union of $K_k^r$ and $(K_{n-k}^r)^c$ by adding all edges with cardinality $r$ connecting $V(K_k^r)$ with $V((K_{n-k}^r)^c)$.

Since $V(K_k^r)$ is a vertex-cut of $H$  and $H-V(K_k^r)=(K_{n-k}^r)^c$, there is no subhypergraph with vertex-connectivity at least $k+1$, and so $\overline{\kappa}(H)\leq k$. If $E(H^c)=\emptyset$, then $H$ is vertex-$k$-maximal by the definition of vertex-$k$-maximal hypergraph. If $E(H^c)\neq\emptyset$, then for any $e\in E(H^c)$, $e$ must be contained in $V((K_{n-k}^r)^c)$ , and so $(H+e)[V(K_k^r)\cup e]$ is isomorphic to $K_{k+r}^r$ and $\kappa((H+e)[V(K_k^r)\cup e])=k+1$. That is $\overline{\kappa}(H+e)\geq k+1$. Thus $H$ is vertex-$k$-maximal.

($ii$) holds by a direct calculation.
$\Box$

\begin{thm}
Let $n,k,r$ be integers such that $k,r\geq2$ and $n\geq k+1$. If  $H$ is vertex-$k$-maximal, then $|E(H)|\geq(^n_r)-(^{n-k}_r)$.
\end{thm}

\noindent{\bf Proof.} We will prove the theorem by induction on $n$. If $n\leq k+r-1$, then by $H$ is vertex-$k$-maximal, we have $H\cong K_n^r$. Thus
$|E(H)|=(^{n}_{r})=(^n_r)-(^{n-k}_r)$ by $n-k\leq r-1$.

Now we assume that $n\geq k+r$, and that the theorem holds for smaller value of $n$. Since $H$ is vertex-$k$-maximal and $n\geq k+r$, we have $H$ is not complete. By Lemma 2.1, $\overline{\kappa}(H)=\kappa(H)=k$, and so $H$ has a separation triple $(S,H_1,H_2)$ with $|S|=k$. Let $n_1=|V(H_1)|$ and $n_2=|V(H_2)|$. Then $n_1,n_2\geq k+1$ and $n=n_1+n_2-k$.

Since $H$ is vertex-$k$-maximal, for any $e\in E((H[S])^c)$, there is a $(k+1)$-vertex-connected subhypergraph $H'$ of $H+e$. By Lemma 2.2, $H'$ is either a subhypergraph of $H_1+e$ or a subhypergraph $H_2+e$. Define

\ \ \ \ \ \ \ \ \ \ \ \ \ \
$E_1=\{e:e\in E((H[S])^c)$ and $ \overline{\kappa}(H_1+e)=k\}$

\ \ \ \ \ \ \ \ \ \ \ \ \ \
$E_2=\{e:e\in E((H[S])^c)$ and $ \overline{\kappa}(H_2+e)=k\}$

\noindent{\bf Claim.} Each of the following holds.

($i$) $E_1\cap E_2=\emptyset$ and $E_1\cup E_2\subseteq E((H[S])^c)$.

($ii$) There is a subset $E_1'\subseteq E_1$ such that $H_1+E_1'$ is vertex-$k$-maximal.

($iii$) There is a subset $E_2'\subseteq E_2$ such that $H_2+E_2'$ is vertex-$k$-maximal.

By the definition, $E_1\cup E_2\subseteq E((H[S])^c)$. Since $H$ is vertex-$k$-maximal, we have $E_1\cap E_2=\emptyset$, and so Claim ($i$) holds.

Assume first that $H_1+E_1$ is complete. If $n_1\leq k+r-1$, then $\overline{\kappa}(H_1+E_1)\leq k$, and so $H_1+E_1$ is vertex-$k$-maximal by the definition of vertex-$k$-maximal hypergraphs. If $n_1\geq k+r$, then by $\overline{\kappa}(H_1)\leq\overline{\kappa}(H)\leq k$ and $\overline{\kappa}(H_1+E_1)\geq k+1$, we can choose a maximum subset $E_1'\subseteq E_1$ such that $\overline{\kappa}(H_1+E_1')\leq k$. It follows by the maximality of $E_1'$ and by the definition of vertex-$k$-maximal hypergraphs that $H_1+E_1'$ is vertex-$k$-maximal. Next, we assume $H_1+E_1$ is not complete. Take an arbitrary edge $e\in E((H_1+E_1)^c)$.  Then $e\in E(H^c)$, and so as $H$ is vertex-$k$-maximal, $H+e$ contains a $(k+1)$-vertex-connected subhypergraph $H'$ with $e\in E(H')$.
If $e\cap(V(H_1)-S)\neq\emptyset$, then by Lemma 2.2,  $H'$ is a subhypergraph of $H_1+e$. If $e\subseteq S$, then as $e\notin E_1$, we can choose $H'$ such that $H'$ is a subhypergraph of $H_1+e$.  That is, $\overline{\kappa}(H_1+E_1+e)\geq k+1$. If $\overline{\kappa}(H_1+E_1)\leq k$, then $H_1+E_1$ is vertex-$k$-maximal.  If $\overline{\kappa}(H_1+E_1)\geq k+1$, then by $\overline{\kappa}(H_1)\leq\overline{\kappa}(H)\leq k$, we can choose a maximum subset $E_1'\subseteq E_1$ such that $\overline{\kappa}(H_1+E_1')\leq k$. It also follows by the maximality of $E_1'$ and by the definition of vertex-$k$-maximal hypergraphs that $H_1+E_1'$ is vertex-$k$-maximal. This verifies Claim ($ii$).
By symmetry, Claim ($iii$) holds. Thus the proof of the Claim is complete.

By Claim ($ii$) and Claim ($iii$), there are $E_1'\subseteq E_1$ and $E_2'\subseteq E_2$ such that $H_1+E_1'$ and $H_2+E_2'$ are vertex-$k$-maximal. Since $n_1,n_2\geq k+1$, by induction assumption, we have $|E(H_1+E_1')|\geq (^{n_1}_r)-(^{n_1-k}_r)$ and $|E(H_2+E_2')|\geq (^{n_2}_r)-(^{n_2-k}_r)$.
By Claim ($i$) and the definition of $(H[S])^c$, we have $|E_1'|+|E_2'|+|E(H[S])|\leq|E_1|+|E_2|+|E(H[S])|\leq|E((H[S])^c)|+|E(H[S])|
=(^k_r)$. Thus

$|E(H)|=|E(H_1)|+|E(H_2)|-|E(H[S])|+|E_H[V(H_1)-S, S, V(H_2)-S]|$

\ \ \ \ \ \ \ \ \
$=|E(H_1+E_1')|-|E_1'|+|E(H_2+E_2')|-|E_2'|-|E(H[S])|+|E_H[V(H_1)-S, S, V(H_2)-S]|$

\ \ \ \ \ \ \ $\geq (^{n_1}_r)-(^{n_1-k}_r)+(^{n_2}_r)-(^{n_2-k}_r)-(^k_r)$

\ \ \ \ \ \ \ \ \ \ $+(^{n}_r)-(^{n_1}_r)-(^{n_2}_r)+(^{k}_r)
-(^{n-k}_r)+(^{n_1-k}_r)+(^{n_2-k}_r)$ (By Lemma 2.3)

\ \ \ \ \ \ \  $=(^n_r)-(^{n-k}_r)$.

This proves Theorem 3.2.
$\Box$

By Lemma 3.1, the lower bound of the sizes of vertex-$k$-maximal hypergraphs given in Theorem 3.2 is best possible. If $r=2$, then a $r$-uniform hypergraph $H$ is just a graph. Thus Theorem 1.2 is a corollary of Theorem 3.2.

\begin{cor} (Xu, Lai and Tian \cite{Xu})
Let $n,k$ be integers with $n\geq k+1\geq3$. If $G$ is a vertex-$k$-maximal graph on $n$ vertices, then $|E(G)|\geq(^n_2)-(^{n-k}_2)= (n-k)k+\frac{k(k-1)}{2}$. Furthermore, this bound is best possible.
\end{cor}

\section{The upper bound of the sizes of vertex-$k$-maximal $r$-uniform hypergraphs}

\noindent{\bf Definition 2.} Let $n,k,r$ be integers such that $k,r\geq2$ and $n\geq 2k$. Assume $n=pk+q$, where $p,q$ are integers and  $1\leq q\leq k$.
We define $H_U(n;k,r)$ to be $((p-1)K_k^r\cup K_q^r) \vee_r (K_k^r)^c$,
where $(p-1)K_k^r$ is the union of $p-1$ complete $r$-uniform hypergraphs on $k$ vertices.

\begin{lem}
Let $n,k,r$ be integers such that $k,r\geq2$ and $n\geq 2k$. If $H= H_U(n;k,r)$, then

($i$) $H$ is vertex-$k$-maximal, and

($ii$) $|E(H)|\leq (^n_r)-(^{n-k}_r)+(\frac{n}{k}-2)(^k_r)$, where the equality holds if $n$ is a multiple of $k$.
\end{lem}

\noindent{\bf Proof.} ($i$) By Definition 2, $H=((p-1)K_k^r\cup K_q^r) \vee_r (K_k^r)^c$. Denote the $p-1$ complete $r$-uniform hypergraphs on $k$ vertices by $K_k^r(1),\cdots,K_k^r(p-1)$. Let $H_0=H[V((K_k^r)^c)]$, $H_p=H[V(K_q^r)]$ and $H_i=H[V(K_k^r(i))]$ for $1\leq i\leq p-1$. Then $H=H_0\vee_r (H_1\cup\cdots\cup H_p)$.

Since $V(H_0)$ is a vertex-cut of size $k$ and every component of $H-V(H_0)$ has at most $k$ vertices. It follows that $H$ contains no $(k+1)$-vertex-connected subhypergraphs, and so $\overline{\kappa}(H)\leq k$. If $E(H^c)=\emptyset$, then $H$ is vertex-$k$-maximal by the definition of vertex-$k$-maximal hypergraphs. Thus we assume $E(H^c)\neq\emptyset$ in the following. Let $e\in E(H^c)$. If $e\subseteq V(H_0)$, then $H'=H[V(H_1)\cup e]$ is isomorphic to $K_{k+r}^r$, and so $\kappa(H')=k+1$. If $e\subseteq V(H_1)\cup\cdots\cup V(H_p)$, let $e$ be exact-($V(H_{i1}),\cdots, V(H_{is})$)-crossing. We will prove that $H''=H[V(H_0)\cup V(H_{i1})\cup\cdots \cup(H_{is})]+e$ is $(k+1)$-vertex-connected. It suffices to prove that $H''-S$ is connected for any $S\subseteq V(H'')$ with $|S|=k$. If $S=V(H_0)$, then, by $e$ is exact-($V(H_{i1}),\cdots, V(H_{is})$)-crossing, $H''-S$ is connected. So assume $V_0'=V(H_0)\setminus S\neq\emptyset$. Let $V_1'=(V(H_{i1})\cup\cdots \cup V(H_{is}))\setminus S$. Then $H''-S$ is isomorphic to $H[V_0']\vee_r H[V_1']$ if $S\cap e\neq\emptyset$;  and $H''-S$ is isomorphic to $H[V_0']\vee_r H[V_1']+e$ if $S\cap e=\emptyset$. Since $V_0', V_1'\neq\emptyset$ and $|V_0'\cup V_1'|\geq r$, we obtain that $H''-S$ is connected. Thus $\overline{\kappa}(H+e)\geq k+1$ for any $e\in E(H^c)$, and so $H$ is vertex-$k$-maximal.

($ii$) By a direct calculation, we have
$|E(H)|\leq (^n_r)-(^{n-k}_r)+(\frac{n}{k}-2)(^k_r)$, where the equality holds if $n$ is a multiple of $k$.
$\Box$

Motivated by Conjecture 1, we propose the following conjecture for vertex-$k$-maximal $r$-uniform hypergraphs.

\begin{conj}
Let $k,r$ be integers with $k,r\geq2$. Then for sufficiently large $n$, every vertex-$k$-maximal $r$-uniform hypergraph $H$ on $n$ vertices satisfies $|E(H)|\leq(^n_r)-(^{n-k}_r)+(\frac{n}{k}-2)(^k_r)$.
\end{conj}

The following theorem confirms Conjecture 2 for the case $k<r$.

\begin{thm}
Let $n,k,r$ be integers such that $k,r\geq2$ and $n\geq 2k$. If $k<r$,
then every vertex-$k$-maximal $r$-uniform hypergraph $H$ on $n$ vertices satisfies $|E(H)|\leq(^n_r)-(^{n-k}_r)+(\frac{n}{k}-2)(^k_r)=(^n_r)-(^{n-k}_r)$.
\end{thm}

\noindent{\bf Proof.}
We will prove the theorem by induction on $n$. If $n\leq k+r-1$, then by $H$ is vertex-$k$-maximal, we have $H\cong K_n^r$. Thus
$|E(H)|=(^{n}_{r})=(^n_r)-(^{n-k}_{r})$ by $n-k\leq r-1$.

Now we assume that $n\geq k+r$, and that the theorem holds for smaller value of $n$. Since $H$ is vertex-$k$-maximal and $n\geq k+r$, we have $H$ is not complete. Let $S$ be a minimum vertex-cut of $H$. By Lemma 2.1, $|S|=k$.
Let $C_1$ be a minimum component of $H-S$ and $C_2=H-(V(C_1)\cup S)$. Assume $H_1=H[V(C_1)\cup S]$ and $H_2=H[V(C_2)\cup S]$.  Since $k<r$, we have $E((H[S])^c)=\emptyset$, and so $H_1$ and $H_2$ are both vertex-$k$-maximal by Lemma 2.2. Let $n_1=|V(H_1)|$ and $n_2=|V(H_2)|$. Then $n=n_1+n_2-k$ and $k+1\leq n_1\leq n_2$. We consider two cases in the following.

\noindent{\bf Case 1.}  $|V(C_1)|=1$.

By $|V(C_1)|=1$, we obtain that $n_2=n-1\geq k+r-1\geq 2k$. Since $H_2$ is vertex-$k$-maximal, by induction assumption, we have
$|E(H_2)|\leq(^{n-1}_r)-(^{n-k-1}_r)$. Thus

$|E(H)|=|E(H_1)|+|E(H_2)|-|E(H[S])|+|E_H[V(H_1)-S, S, V(H_2)-S]|$

\ \ \ \ \ \ \ \ \ \
$\leq (^k_{r-1})+(^{n-1}_r)-(^{n-k-1}_r)
+(^{n-1}_{r-1})-(^{k}_{r-1})-(^{n-k-1}_{r-1})$

\ \ \ \ \ \ \ \ \ \
$= (^n_r)-(^{n-k}_r)$.

\noindent{\bf Case 2.} $|V(C_1)|\geq 2$.

By $|V(C_1)|\geq 2$, we obtain that $C_1$ contains edges, and so $|V(C_1)|\geq r$. Thus $n_2\geq n_1\geq k+r\geq 2k+1$. Since both $H_1$ and $H_2$ are vertex-$k$-maximal, by induction assumption, we have
$|E(H_i)|\leq(^{n_i}_r)-(^{n_i-k}_r)$ for $i=1,2$. Thus

$|E(H)|=|E(H_1)|+|E(H_2)|-|E(H[S])|+|E_H[V(H_1)-S, S, V(H_2)-S]|$

\ \ \ \ \ \ \ \ \ \
$\leq (^{n_1}_r)-(^{n_1-k}_r)+(^{n_2}_r)-(^{n_2-k}_r)$

\ \ \ \ \ \ \ \ \ \ \ \
$+(^{n}_r)-(^{n_1}_r)-(^{n_2}_r)+(^{k}_r)
-(^{n-k}_r)+(^{n_1-k}_r)+(^{n_2-k}_r)$ (By Lemma 2.3)

\ \ \ \ \ \ \ \ \ \
$=(^n_r)-(^{n-k}_r)$.

This completes the proof.
$\Box$

Combining Theorem 3.2 with Theorem 4.2, we have the following corollary.

\begin{cor}
Let $n,k,r$ be integers such that $k,r\geq2$ and $n\geq 2k$. If $k<r$,
then every vertex-$k$-maximal $r$-uniform hypergraph $H$ on $n$ vertices satisfies $|E(H)|=(^n_r)-(^{n-k}_r)$.
\end{cor}


\vspace{1cm}
\begin{thebibliography}{5}

\bibitem{Anderson} J. Anderson, H.-J. Lai, X. Lin, M. Xu, On $k$-maximal strength digraphs, J. Graph Theory 84 (2017) 17-25.

\bibitem{Bahmanian} M. A. Bahmanian, M. $\breve{S}$ajna, Connection and separation in hypergraphs, Theory and Applications of Graphs 2(2) (2015) 0-24.

\bibitem{Bernshteyn} A. Bernshteyn, A. Kostochka, On the number of edges in a graph with no ($k+1$)-connected subgraphs, Discrete Mathematics 339 (2016) 682-688.

\bibitem{Bondy} J. A. Bondy, U. S. R. Murty, Graph Theory, Graduate Texts in Mathematics 244, Springer, Berlin, 2008.

\bibitem{Chekuri} C. Chekuri, C. Xu, Computing minimum cuts in hypergraphs, In Proceedings of the Twenty-Eighth Annual ACM-SIAM Symposium on Discrete Algorithms (Barcelona, 2017), SIAM (2017) 1085-1100.

\bibitem{Dewar} M. Dewar, D. Pike, J. Proos, Connectivity in Hypergraphs,  arXiv:1611.07087v3.

\bibitem{Lai}  H.-J. Lai, The size of strength-maximal graphs, J. Graph Theory 14 (1990) 187-197.

\bibitem{Lin} X. Lin, S. Fan, H.-J. Lai, M. Xu, On the lower bound of $k$-maximal digraphs, Discrete Math. 339 (2016) 2500-2510.

\bibitem{Mader71} W. Mader, Minimale $n$-fach kantenzusammenhngende graphen, Math. Ann. 191 (1971) 21-28.

\bibitem{Mader72} W. Mader, Existenz n-fach zusammenh$\ddot{a}$ngender Teilgraphen in Graphen gen$\ddot{u}$gend gro$\beta$en Kantendichte,
Abh. Math. Sem. Univ. Hamburg 37 (1972) 86-97.

\bibitem{Mader79} W. Mader, Connectivity and edge-connectivity in finite graphs, in: B. Bollob$\acute{a}$s (Ed.), Surveys
in Combinatorics, Cambridge University Press, London, 1979, pp. 66-95.

\bibitem{Matula69} D. W. Matula, The cohesive strength of graphs, The Many Facets of Graph Theory, Lecture Notes in Mathematics, No. 110, G. Chartrand and S. F. Kapoor, eds., Springer-Verlag, Berlin, 1969, pp. 215-221.

\bibitem{Matula72} D. Matula, $K$-components, clusters, and slicings in graphs, SIAM J. Appl. Math. 22 (1972) 459-480.

\bibitem{Matula78} D. W. Matula, Subgraph connectivity numbers of a graph, Theory and Applications of Graphs, pp. 371-383 (Springer-Verlag, Berlin 1978).

\bibitem{Matula83} D. W. Matula, Ramsey theory for graph connectivity, J. Graph Theory 7 (1983) 95-105.

\bibitem{Tian} Y. Z. Tian, L. Q. Xu, H.-J. Lai, J. X. Meng, On the sizes of $k$-edge-maximal $r$-uniform hypergraphs, arXiv:1802.08843v3.

\bibitem{Xu} L. Q. Xu, H.-J. Lai, Y. Z. Tian, On the sizes of vertex-$k$-maximal graphs, submitted.

\bibitem{Yuster} R. Yuster, A note on graphs without $k$-connected subgraphs, Ars Combin. 67 (2003) 231-235.


\end{thebibliography}
\end{document}